\documentclass[11pt,leqno]{article}
\usepackage{amsmath,amsthm,amsfonts,amssymb,amscd}

{\newtheorem{Theorem}{Theorem}
\newtheorem{Corollary}{Corollary}
\newtheorem{Proposition}{Proposition}
\newtheorem{Lemma}{Lemma}
\newtheorem{Claim}{Claim}
\theoremstyle{Definition}
\newtheorem{Definition}{Definition}
\newtheorem{Example}{Example}
\newtheorem{Question}{Question}
\newtheorem{Problem}{Problem}

\theoremstyle{Remark}
\newtheorem{Remark}{Remark}

\def\leaderfill{\leaders\hbox to .8em{\hss .\hss}\hfill}
\def\_#1{{\lower 0.7ex\hbox{}}_{#1}}
\def\esima{${}^{\text{\b a}}$}
\def\esimo{${}^{\text{\b o}}$}
\font\bigbf=cmbx10 scaled \magstep2

\def\C{{\mathcal{C}}}
\def\L{{\mathcal{L}}}
\def\G{{\mathcal{G}}}
\def\fa{{\mathcal{F}}}
\def\bj{{\mathbb{J}}}
\def\ba{{\mathbb{A}}}
\def\O{{\mathcal{O}}}
\def\eR{{\mathcal{R}}}
\def\M{{\mathcal{M}}}
\def\D{{\mathcal{D}}}
\def\U{{\mathcal{U}}}
\def\B{{\mathcal{B}}}
\def\E{{\mathcal{E}}}

\def\po{{\partial}}
\def\ro{{\rho}}
\def\te{{\theta}}
\def\Te{{\Theta}}
\def\Om{{\Omega}}
\def\vr{{\varphi}}
\def\ga{{\gamma}}
\def\Ga{{\Gamma}}
\def\la{{\lambda}}
\def\La{{\Lambda}}
\def\ov{\overline}
\def\al{{\alpha}}
\def\epsilon{{\varepsilon}}
\def\lg{{\langle}}
\def\rg{{\rangle}}
\def\lv{{\left\epsilonrt}}
\def\rv{{\right\epsilonrt}}

\def\A{{\mathbb{A}}}
\def\bh{{\Bbb{H}}}
\def\bp{{\Bbb{P}}}
\def\ee{{\Bbb{E}}}
\def\re{{\Bbb{R}}}
\def\bz{{\Bbb{Z}}}
\def\bq{{\Bbb{Q}}}
\def\bd{{\Bbb{D}}}
\def\bc{{\Bbb{C}}}
\def\bn{{\Bbb{N}}}
\def\be{{\beta}}
\def\Det{\operatorname{{Det}}}
\def\dim{\operatorname{{dim}}}

\def\grad{\operatorname{{grad}}}
\def\cod{\operatorname{{cod}}}
\def\Ind{\operatorname{{Ind}}}
\def\SL{\operatorname{{SL}}}
\def\Res{\operatorname{{Res}}}
\def\Fol{\operatorname{{Fol}}}
\def\tr{\operatorname{{tr}}}
\def\dim{\operatorname{{dim}}}
\def\Aut{\operatorname{{Aut}}}
\def\SL{\operatorname{{SL}}}
\def\Aff{\operatorname{{Aff}}}
\def\Hol{\operatorname{{Hol}}}
\def\loc{\operatorname{{loc}}}
\def\Ker{\operatorname{{Ker}}}
\def\mI{\operatorname{{Im}}}
\def\Dom{\operatorname{{Dom}}}
\def\Id{\operatorname{{Id}}}
\def\Tni{\operatorname{{Int}}}
\def\supp{\operatorname{{supp}}}
\def\Diff{\operatorname{{Dif}}}
\def\sing{\operatorname{{sing}}}

\title{On the classification of non-integrable complex distributions}

\author{Toshikazu Ito  \& Bruno Sc\'ardua}
\date{}
\begin{document}
\maketitle


\def\leaderfill{\leaders\hbox to .8em{\hss .\hss}\hfill}
\def\_#1{{\lower 0.7ex\hbox{}}_{#1}}
\def\esima{${}^{\text{\b a}}$}
\def\esimo{${}^{\text{\b o}}$}
\font\bigbf=cmbx10 scaled \magstep2

\def\A{{\mathcal{A}}}
\def\sa{{\mathcal{S}}}
\def\P{{\mathcal{P}}}
\def\C{{\mathcal{C}}}
\def\L{{\mathcal{L}}}
\def\G{{\mathcal{G}}}
\def\fa{{\mathcal{F}}}
\def\O{{\mathcal{O}}}
\def\eR{{\mathcal{R}}}
\def\M{{\mathcal{M}}}
\def\D{{\mathcal{D}}}
\def\U{{\mathcal{U}}}
\def\B{{\mathcal{B}}}
\def\E{{\mathcal{E}}}

\def\po{{\partial}}
\def\ro{{\rho}}
\def\te{{\theta}}
\def\Te{{\Theta}}
\def\Om{{\Omega}}
\def\vr{{\varphi}}
\def\ga{{\gamma}}
\def\Ga{{\Gamma}}
\def\la{{\lambda}}
\def\La{{\Lambda}}
\def\ov{\overline}
\def\al{{\alpha}}
\def\epsilon{{\varepsilon}}
\def\lg{{\langle}}
\def\rg{{\rangle}}
\def\lv{{\left\epsilonrt}}
\def\rv{{\right\epsilonrt}}
\def\be{{\beta}}

\def\bh{{\Bbb{H}}}
\def\bp{{\Bbb{P}}}
\def\ee{{\Bbb{E}}}
\def\re{{\Bbb{R}}}
\def\bz{{\Bbb{Z}}}
\def\bq{{\Bbb{Q}}}
\def\bd{{\Bbb{D}}}
\def\bc{{\Bbb{C}}}
\def\bn{{\Bbb{N}}}

\def\SL{\operatorname{{SL}}}
\def\Res{\operatorname{{Res}}}
\def\Fol{\operatorname{{Fol}}}
\def\tr{\operatorname{{tr}}}
\def\dim{\operatorname{{dim}}}
\def\Aut{\operatorname{{Aut}}}
\def\GL{\operatorname{{GL}}}
\def\Aff{\operatorname{{Aff}}}
\def\Hol{\operatorname{{Hol}}}
\def\loc{\operatorname{{loc}}}
\def\Ker{\operatorname{{Ker}}}
\def\mI{\operatorname{{Im}}}
\def\Dom{\operatorname{{Dom}}}
\def\Id{\operatorname{{Id}}}
\def\Tni{\operatorname{{Int}}}
\def\supp{\operatorname{{supp}}}
\def\Diff{\operatorname{{Diff}}}
\def\sing{\operatorname{{sing}}}
\def\Sing{\operatorname{{sing}}}
\def\codim{\operatorname{{codim}}}
\def\grad{\operatorname{{grad}}}
\def\Ind{\operatorname{{Ind}}}
\def\deg{\operatorname{{deg}}}
\def\rank{\operatorname{{rank}}}

\section{Introduction}
\label{section:Introduction} Given a point $p\in {\mathbb C}^n$
and a real number $r>0$ we denote by $B^{2n}(p,r)$ the open ball
of radius $r$ centered at $p$ in ${\mathbb C}^{n}$. The
corresponding closed ball is denoted by $B^{2n}[p,r]$ and its
boundary sphere by $S^{2n-1}(p,r)= \partial B^{2n}[p,r]$.  We also
write $B^{2n}(1)=B^{2n}(0,1)$, $B^{2n}[1]=B^{2n}[0,1]$ and
$S^{2n-1}(1)=\partial B^{2n}[1]$. Let $\Omega$ be a germ of
holomorphic one-form with an isolated singularity  at the origin
$0\in \bc ^n, \, n\ge 3$. We address the problem of analytical
classification of $\Omega$ in the {\it non-integrable} case.
Motivated by the geometrical-analytical classification of
singularities in dimension $2$ we consider the case where the
kernel of $\Omega$ generates a germ of distribution $\Ker(\Omega)$
transverse to small spheres $S^{2n-1}(0,\epsilon)$. This is one,
though not the only, central motivation for this work. The problem
of existence of integral manifolds for germs of singularities of
integrable one-forms is an ancient problem already considered in
the work of Briot-Bouquet. The existence results in
\cite{Camacho-Sad} (for dimension $n=2$) and in
\cite{Cano-Cerveau} (for the non-dicritical case in dimension
$n=3$)  motivate the very basic question below:

\begin{Question}
\label{Question:first}  Is there a non-integrable germ of
holomorphic one-form $\Omega$ with an isolated singularity  at the
origin $0 \in \bc^{n}$ such that $\Ker(\Omega)$ is transverse to
the spheres $S^{2n-1}(0,\epsilon)$, \, for $\epsilon > 0$ small
enough and $\Ker(\Omega)$ admits no integral manifold through the
origin?
\end{Question}

Theorem~\ref{Theorem:nointegral manifold} gives a positive answer
to this question.  Other motivations are related to our previous
work  in \cite{Ito-Scardua1} and  \cite{Ito-Scardua2} where we
study the obstructions to the integrability of $\Omega$.
 Our first main result reads as:

\begin{Theorem}
\label{Theorem:main}  Let $n\ge 3$ and $\Omega$ be a holomorphic
one-form defined in a neighborhood of $B^{2n}[1]$  and such that
$\sing(\Omega) \cap S^{2n-1}(1)=\emptyset$. If  there exists a
holomorphic vector field $\xi$ in a neighborhood of $B^{2n}[1]$,
transverse to $S^{2n-1}(1)$, and such that $\Omega \cdot\xi= 0$,
then   $\Omega$ is not integrable.
\end{Theorem}

 Let us give examples of
 distributions as in  Theorem~\ref{Theorem:main}.
 Denote by $\mathcal A(2m)$ the set of all $2m\times2m$
skew-symmetric complex matrices and by $\ba(2m)$ the subset of
nonsingular elements in $\mathcal A(2m)$.
 In \cite{Ito-Scardua2} it is observed that if $A =
(a_{ij})_{i,j=1}^{2m}$  belongs to $\mathbb A(2m)$  then, for $m
\ge 2$,\,the one-form\, $\Omega_A = \sum\limits_{i,j=1}^{2m}
a_{ij}\, z_i\,dz_j$ defines a non-integrable holomorphic (linear)
distribution transverse to the spheres $S^{4m-1}(0,r) \subset \bc
^{2m}$, \, $r > 0$. Such a one-form will be called {\it linear}. A
particular case is the one-form $\Omega_{\bj(2m)} =
\sum\limits_{j=1}^m (z_{2j-1}\,dz_{2j} - z_{2_j}\,dz_{2j-1})$. One
may ask for non-linear examples. Given
$\ell=(\ell_1,...,\ell_m)\in {\mathbb N}^m$, we introduce the
corresponding {\it non-integrable Poincar\'e-Dulac normal form} as
$\Omega_{(\ell)} = \sum\limits_{j=1}^m [z_{2j-1}\,dz_{2j} -
(\ell_j z_{2j} + z_{2j-1}^{\ell_j})dz_{2j-1}]$ in coordinates
$(z_1,z_2,\dots,z_{2m}) \in \bc^{2m}$. We prove that also
$\Omega_{(\ell)}$ is  not integrable for $m \ge 2$, singular only
at the origin,  and $\Ker(\Omega_{(\ell)})$ is transverse to
$S^{4m-1}(0,r)$, \, $\forall\, r>0$ small enough (see
Example~\ref{Example:distributions}). The one-forms $\Omega_A,
\Omega_{\mathbb J(2m)}$ and $\Omega_{(\ell)}$ are our basic models
in the classification we pursue (see \S 4).

Let $\Ker(\Omega)$ be a codimension one holomorphic distribution
on a complex manifold $V^n$. Let  $p \in V^n$ be a {\it
singularity\/} of $\Ker(\Omega)$, that is, of $\Omega$. A germ of
codimension one analytic subset $\La_p$ at $p$ is  an {\it
integral manifold\/} of $\Ker(\Omega)$ through $p$ if any vector
$v_q\in T_q(V)$ which is tangent to $\Lambda_p$ at a point $q$
belongs to $\Ker(\Omega)(q)$. This means that if $\La$ is any
representative of $\La_p$ in a neighborhood $U$ of $p$ in $V^n$
and $\Lambda^*$ denotes the smooth part of $\Lambda$ then the
tangent bundle $T\La^*$ is a sub-bundle of
$\Ker(\Omega)\big|_{\Lambda^*}$ (see
Definition~\ref{Definition:invariant}). We shall always assume
$\La$ and $\La_p$ to be irreducible, nevertheless we do not
require that $\La\backslash\La^* = \sing(\La)$ is contained in
$\Sing(\Ker(\Omega))$. Regarding the existence of integral
manifolds for non-integrable distributions we have:

\begin{Theorem}
\label{Theorem:nointegral manifold} Let $m\ge 2$. Given $A\in
\mathbb A(2m)$ and $\ell\in {\mathbb N}^m$ the distributions
$\Ker(\Omega_{\mathbb J(2m)}), \, \Ker(\Omega_A)$ and
$\Ker(\Omega_{(\ell)})$ admit no integral manifold through the
origin.
\end{Theorem}

In the course of the proof of Theorem~\ref{Theorem:nointegral
manifold} we obtain the following Darboux's theorem type for (not
necessarily integrable) polynomial distributions. This is actually
a non-integrable version of the more precise Theorem 3.3 on page
102 of \cite{Jouanolou}:

\begin{Proposition}
\label{Proposition:Darboux}  Let $\Omega$ be a {\rm(}not
necessarily integrable{\rm)} polynomial one-form on $\bc^n$, \, $n
\ge 2$ and assume that $\cod \sing(\Omega)\ge 2$. If
$\Ker(\Omega)$ has infinitely many algebraic invariant
hypersurfaces then $\Omega$ is integrable. Indeed $\Omega = PdQ -
QdP$ for some polynomials $P,Q$ with no common factors and, in
particular, the leaves of the foliation $\fa_\Omega$ defined by
$\Omega$ are contained in the algebraic subvarieties $\{\la P -
\mu Q=0\}$ where $(\la,\mu) \in \bc^2-\{(0,0)\}$.
\end{Proposition}

As already mentioned, the examples  $\Omega_A, \,\Omega_{\bj(2m)}$
and $\Omega_{(\ell)}$ above constructed motivate the problem of
analytical classification of germs of non-integrable one-forms
defining distributions transverse to small spheres (see
Questions~\ref{Question:final} and \ref{Question:classify} in
\S\ref{section:classification}). In this direction we prove:

\begin{Theorem}
\label{Theorem:isotopy}   Let $\Om$ be a holomorphic one-form in a
neighborhood $U$ of the closed ball  $B^{4m}[1]\subset \bc^{2m}$
and such that  {\rm(1)} $\Om\cdot\vec R = 0$, where $\vec R$ is
the radial vector field in $\bc^{2m}$  and {\rm(2)} $\Sing(\Om)
\cap S^{4m-1}(1) = \emptyset$. Then $\Ker(\Om)$ is homotopic to
the linear distribution $\Ker(\Omega_{\bj(2m)})$ by distributions
$\Ker(\Om_s)$, $0 \le s \le 1$, such that $\Omega_0=\Omega$ and
$\Omega_1=\Omega_{\mathbb J_ {2m}}$, where $\Om_s$ is holomorphic
and satisfies $(1)$ and $(2)$ above.
\end{Theorem}

\noindent{\bf Acknowledgement}. We want to thank the anonymous
referee for  many valuable remarks, suggestions and for pointing
out gaps and unclear points in the first version. Specially for
suggesting clearer arguments and in particular for a drastic
simplification in the proof of Theorem~\ref{Theorem:nointegral
manifold}. Also mainly due to the referee is the discussion in
Remark~\ref{Remark:final}.

\section{Invariant manifolds}
\label{section:integral manifolds} In this section we discuss
Question~\ref{Question:first} and prove
Theorem~\ref{Theorem:nointegral manifold} and
Proposition~\ref{Proposition:Darboux}. First we prove the examples
mentioned in the introduction.

\begin{Example}
\label{Example:distributions}
 {\rm Let $A = (a_{ij})\in \mathbb A(2m)$ and
 $\Omega_A := \sum\limits_{i,j=1}^{2m} a_{ij}z_idz_j$ the
corresponding linear one-form in $\bc^{2m}$. Then $\Sing(\Omega_A)
= \{0\} \subset \bc^{2m}$ and $\Omega_A \cdot \vec R = 0$ for the
radial vector field $\vec R = \sum\limits_{j=1}^{2m}
z_j\,\dfrac{\po}{\po z_j}\cdot$ This implies that $\Ker(\Omega_A)$
is transverse to every sphere $S^{4m-1}(0,r)$, $r > 0$. The
non-integrability of $\Omega_A$ and $\Omega_{(\ell)}$ is a
straightforward computation (cf.\cite{Ito-Scardua1}).  For the
transversality of the non-integrable Poincar\'e-Dulac  normal form
$ \Omega_{(\ell)}$ with small spheres $S^{4m-1}(0,\epsilon)$  we
observe that if $\xi_{(\ell)} = \sum\limits_{j=1}^m
\big[(\ell_jz_{2j} + z_{2j-1}^{\ell_j})\, \frac{\po}{\po z_{2j}} +
z_{2j-1}\, \frac{\po}{\po z_{2j-1}}\big]$
 then $\Omega_{(\ell)} \cdot \xi_{(\ell)}=0$ and, as it is well-known,
 the vector field $\xi_{(\ell)}$ is transverse to the
spheres $S^{4m-1}(0,\epsilon)$ if $\epsilon>0$ is small enough.}
\end{Example}

 Now we shall prove Theorem~\ref{Theorem:nointegral manifold}.
 Let us first precise some notions involved.

\begin{Definition}
\label{Definition:invariant} {\rm Let $\Omega$ be a holomorphic
one-form on a complex manifold $V$. A {\it smooth} complex
immersed  submanifold $\La^* \subset V$ is an {\it integral
manifold\/} of $\Omega$ if  $T\La^* \subset
\Ker(\Omega)|\_{\La^*}$\,. In other words, any tangent vector to
$\La^*$ belongs to $\Ker(\Omega)$. If $\La \subset V$ is a
possibly singular complex analytic submanifold we say that $\La$
is an {\it integral manifold\/} of $\Omega$ if its regular part
$\La^*=\Lambda\setminus \sing(\Lambda)$ is an integral manifold of
$\Omega$. }
\end{Definition}

The following lemma is found in the algebraic setting in
\cite{Jouanolou}, Section 3.1, page 99.

\begin{Lemma}
\label{Lemma:Darboux} Let $\Ker(\Omega)$ be given by a holomorphic
one-form $\Omega$ in $V$ with $\cod \sing(\Omega) \ge 2$ and let
$\La \subset V$ be a codimension one analytic subset given by a
reduced equation $\La : \{f=0\}$ for some holomorphic $f\colon V
\to \bc$. The following conditions are equivalent:

\noindent {\rm (1)} $\La$ is $\Ker(\Omega)$-invariant. \, \, {\rm
(2)} $\Omega \wedge \frac{df}{f}$ is a holomorphic $2$-form on
$V$.

\end{Lemma}

\noindent{\bf Proof}. Since all objects involved are analytic we
may consider  the local case also at a generic (and therefore
non-singular) point $p \in \La^*$. In suitable local coordinates
$(z_1,\dots,z_n) = (z_1,\dots,z_{n-1},f)$ we have $p=0$ and $\La$
given by $\{f=z_n=0\}$. Also we may write $\Omega =
\sum\limits_{j=1}^n a_j\,dz_j$\,.  Suppose $\La$ is $\Ker(\Omega)$
invariant. Then, since $\La$ is given by $\{z_n=0\}$ we have
$\Omega \cdot \frac{\po}{\po z_j}\bigg|_{\{z_n=0\}} = 0$ for all
$j \in \{1,\dots,n-1\}$. In other words $z_n$ divides $a_j$ in
$\bc\{z_1,\dots,z_n\}$ for every $j \in \{1,\dots,n-1\}$ and
therefore $\Omega = \sum\limits_{j=1}^{n-1} z_n\, \tilde a_j\,
dz_j + a_n\, dz_n$ for some holomorphic $\tilde a_1,\dots,\tilde
a_{n-1} \in \bc\{z_1,\dots,z_n\}$. Thus $\Omega \wedge
\frac{dz_n}{z_n} = \sum\limits_{j=1}^{n-1} \tilde a_j\, dz_j
\wedge dz_n$ is holomorphic. Conversely, if  $\Omega \wedge
\frac{dz_n}{z_n}$ is holomorphic then $\Omega =
\sum\limits_{j=1}^{n-1} z_n\,\tilde a_j\,dz_j + a_n\,dz_n$ as
above therefore $\Omega \cdot\frac{\po}{\po z_j}$ vanishes on
$\{z_n=0\}$ for every $j=1,\dots,n-1$. Since $\{\frac{\po}{\po
z_1},\dots, \frac{\po}{\po z_{n-1}}\}$ generate $T\La$ in a
neighborhood of 0 we obtain that $T\La \subset \Ker(\Omega)$ in a
neighborhood of $p$. Therefore $\La$ is $\Ker(\Omega)$-invariant.
\qed

\vglue.1in

\noindent{\bf Proof of Proposition~\ref{Proposition:Darboux}}.
Once we have Lemma~\ref{Lemma:Darboux}, the proof is essentially
the same given in \cite{Jouanolou}, thus we will omit it and refer
to \cite{Jouanolou}.  \qed \vglue.1in

\subsection{Proof of
Theorem~\ref{Theorem:nointegral manifold}}
 Let $\Omega_{\bj(2m)}=\sum\limits_{j=1}^m (z_{2j-1}dz_{2j} - z_{2j}
dz_{2j-1})$ be given in $\bc^{2m}$. Then $d\Omega_{\bj(2m)}=2
\sum\limits_{j=1}^{m}dz_{2j-1}\wedge dz_{2j}$ is a non-degenerate
$2$-form on $\bc^{2m}$ and $(d\Omega_{\bj(2m)})^m$ is a non-zero
$2m$-form. Let $L$ be a complex manifold and $\zeta \colon L \to
\bc^{2m}$ a smooth embedding.

\begin{Lemma}
\label{Lemma:nointegral manifold} The following statements are
equivalent:

\noindent{\rm (1)} $T(\zeta(L))\subset
\Ker(\Omega_{\bj(2m)})\big|_{\zeta(L)}.$

\noindent{\rm (2)} $\zeta^*(\Omega_{\bj(2m)})=0$.
\end{Lemma}

\noindent{\bf Proof}. To prove that (1) implies (2) we take any
vector $v\in T(L)$, then $\zeta_*(v)\in T(\zeta(L))\subset
\Ker(\Omega_{\bj(2m)})\big|_{\zeta(L)}$ means
$\Omega_{\bj(2m)}(\zeta_*v)=0$. Therefore we get
$\zeta^*(\Omega_{\bj(2m)})(v)=\Omega_{\bj(2m)}(\zeta_*v)=0$. To
prove that (2) implies (1) take any $\tilde v = \zeta_* v \in
T(\zeta(L)), \, v \in T(L)$. Then $\zeta^*(\Omega_{\bj(2m)})=0$
means that $\Omega(\tilde
v)=\Omega_{\bj(2m)}(\zeta_*v)=\zeta^*(\Omega_{\bj(2m)})(v)=0$.
Thus $\tilde v$ belongs to
$\Ker(\Omega_{\bj(2m)})\big|_{\zeta(L)}$. \qed

\begin{Proposition}
\label{Proposition:dimensionless} Under the above notations,
assume that $T(\zeta(L))\subset
\Ker(\Omega_{\bj(2m)})\big|_{\zeta(L)}$, then the complex
dimension of $L$ is less than $m$ or equal to $m$, i.e., $\dim
_{\bc} L \leq m$.
\end{Proposition}

\noindent{\bf Proof}. First we observe that
$\zeta^*(d\Omega_{\bj(2m)})= d(\zeta^*\Omega_{\bj(2m)})=0$ on $L$.
Assume that $\dim_\bc L \geq m+1$. Take $m+1$ linearly independent
vectors $v_1,...,v_{m+1}$ in $T(L)$. Moreover, we take $m-1$
vectors $u_1,...,u_{m-1}$ in $T\bc^{2m}$ such that $\zeta_*
v_1,...,\zeta_* v_{m+1}, u_1,...,u_{m-1}$ are linearly independent
in $T\bc^{2m}$. Then we get
$(d\Omega_{\bj(2m)})^m(\zeta_*v_1,...,\zeta_*v_{m+1},u_1,...,u_{m-1})=0$
because $\zeta^*(d\Omega_{\bj(2m)})=0$. It is contradictory with
the fact that $(d\Omega_{\bj(2m)})^m$ is a non-zero $2m$-form.
\qed

\begin{Corollary}
\label{Corollary:nointegral manifold} Let $m\ge 2$. The one-form
$\Omega_{\bj(2m)}$ has no integral manifold through the origin.
\end{Corollary}

\noindent{\bf Proof}. Assume that $\Lambda$ is an integral
manifold of $\Omega_{\bj(2m)}$. Then $\Lambda^*= \Lambda -\{0\}$
is a complex hypersurface, i.e., $\dim_\bc \Lambda ^*=2m-1$. By
Proposition~\ref{Proposition:dimensionless}, we get $\dim_\bc
\Lambda^* \leq m$. We have a contradiction with the hypothesis
$m\ge 2$.

\vglue.2in This proves Theorem~\ref{Theorem:nointegral manifold}
for $\Omega_{\mathbb J(2m)}$.
 By the same argument,
$\Omega_A$ and $\Omega_{(\ell)}$ have no integral manifold through
the origin of $\bc^{2m}$. \qed

\vglue.1in Regarding the smoothness of invariant hypersurfaces of
holomorphic
 distributions with an isolated singularity we have:

\begin{Proposition}
\label{Proposition:smooth}
 Let $\Omega$ be a holomorphic one-form in a
complex manifold $M$ and $p\in M$ and $\Lambda\subset M$ an
integral manifold of $\Omega$ point. Given a point $p\in M$  where
$\Omega$ is nonzero then $\Lambda$  is smooth at the point $p$.
\end{Proposition}

\noindent{\bf Proof}. In suitable local coordinates
$(z_1,...,z_n)$ centered at $p$ the Pfaffian equation $\Omega=0$
is equivalent to an equation $dz_n=\sum\limits_{j=1}^{n-1} g_j
dz_j$ for some holomorphic functions $g_j$. Any integral manifold
then writes as $z_n=f(z_1,...,z_{n-1})$ for a holomorphic function
$f$ satisfying $\po f (z_1,...,z_{n-1})/ \po
z_j=g_j(z_1,...,z_{n-1},f(z_1,...,z_{n-1}))$. This implies that
$f$ cannot develop any singularities. \qed

\section{Proof of Theorem~\ref{Theorem:main}}
\label{section:proofmain} By hypothesis,  $\xi$ defines  a
one-dimensional holomorphic foliation  in a neighborhood of the
closed ball $B^{2n}[1]$, transverse to the sphere $S^{2n-1}(1)$.
 According to \cite{Ito}   $\xi$ has a unique
singularity $p$ in the open ball $B^{2n}(1)$ and this singularity
is in the Poincar\'e domain: the  vector field $\xi$ has a
non-singular linear part at $p$ having eigenvalues
$\la_1,...,\lambda_n$ such that the origin does not belong to the
convex hull of the set $\{\la_1,...,\la_n\}$ in $\re^2$. On the
other hand we have  the following result from \cite{Ito-Scardua2}:

\begin{Theorem} \label{Theorem 8.1}  Let $\omega$ be a holomorphic one-form in
an open subset $U \subset \bc^n$, $n \ge 2$. Suppose that the
distribution $\Ker(\omega)$ is transverse to a sphere
$S^{2n-1}(p,r) \subset U$, with $B^{2n}[p,r] \subset U$, $p \in
U$, $r>0$, then $n$ is even, $\omega$ has a single singular point
in the ball $B^{2n}(p,r)$ and this is a simple singularity in the
following sense:  if we write $\omega=\sum\limits_{j=1}^n f_j
dz_j$ in local coordinates centered at the singularity then the
matrix $(\partial f_j /\partial z_k)_{j,k=1}^n$ is non-singular at
the singularity.
\end{Theorem}

 By the above result  $\Omega$ has a unique singularity $q$ in $B^{2n}(1)$
 and this is a simple singularity.

\begin{Claim}
We have $p=q$.
\end{Claim}
\noindent{\bf Proof}.
 Let $\Omega =
\sum\limits_{j=1}^n f_j\,dz_j$  in a neighborhood $U$ of
$B^{2n}[1]$ on $\bc^n$. Write $\xi = \sum\limits_{i=1}^n
A_i\,\frac{\po}{\po z_i}$ then $0 = \Omega\cdot\xi =
\sum\limits_{j=1}^n f_jA_j$ and therefore $0 = \frac{\po}{\po z_k}
(\sum_{j=1}^n f_jA_j) = \sum_{j=1}^n (\frac{\po f_j}{\po z_k}\,
A_j + f_j\, \frac{\po A_j}{\po z_k}). $ Thus $0 =
\sum\limits_{j=1}^n \frac{\po f_j}{\po z_k}(p)\cdot A_j(p)$, \,
$\forall\,k$ and since the matrix $\big(\frac{\po f_j}{\po
z_k}(p)\big)_{j,k=1}^n$ is non-singular we have $A_j(p)=0$, \,
$j=1,\dots,n$. By the uniqueness of the singularity of $\xi$ we
get $p=q$. This proves the claim.

\vglue.1in Let us finish the proof.  Suppose by contradiction that
$\Omega$ is integrable. By a theorem of Malgrange
\cite{[Malgrange]}, since $\codim \sing(\Omega)=n \ge 3$ at $p$,
the one-form $\Omega$ admits a holomorphic (Morse-type) first
integral in a neighborhood of $p$, say $f\colon (W,p) \to
(\bc,0)$. Then $\Omega\cdot \xi=0$ implies that $\xi(f)=0$.
Because the germ of $\xi$ at $0$ is in the Poincar\'e domain, this
implies that $f$ is constant in a neighborhood of $0$,
contradiction with the fact that $f$ is of Morse type at $0$. This
ends the proof of Theorem~\ref{Theorem:main}. \qed

\section{On the analytical classification}
\label{section:classification}

The problem of  analytic classification of singularities of
 holomorphic one-forms in dimension two is a very well-developed
 topic.  Recently (cf. \cite{Cerveau-Mozo}) the analytic classification
 was obtained for germs of reduced {\it integrable} one-forms at
 the origin $0\in \bc^3$. As far as
 we know, nothing is found regarding the non-integrable case.
In the follow-up,  "to classify"   means to give a description in
terms of objects which are completely understood.
 Obviously  the class of germs of singular (non-integrable) holomorphic
 one-forms is too wide in order to be classified at a first
 moment and we also lack of geometric ingredients. This remark
 is one of the motivations for our  approach in this section. Other motivation is
 given by the well-known results for holomorphic foliations with singularities
 in dimension two collected in the
 following omnibus theorem:

\begin{Theorem}[\cite{[Brjuno]},
\cite{[Camacho-Kuiper-Palis]}, \cite{[Dulac]}, \cite{Ito}]
Given a germ of singular holomorphic one-form $\Omega$ at $0 \in
\bc^2$ the following conditions are equivalent:
\begin{enumerate}
\item[{\rm(1)}] $\Omega$ is in the {\it Poincar\'e domain}:
$\Omega = Ady - Bdx$ where the vector field $X = A\,\frac{\po}{\po
x} + B\,\frac{\po}{\po y}$ has a singularity at $(0,0)$ whose
linear part $DX(0,0)$  has non-zero eigenvalues $\la$, $\mu$ with
quotient $\la/\mu \in \bc\backslash\re_-$\,. \item[{\rm(2)}]
$\Ker(\Omega)$ is transverse to some {\rm(}and therefore to
every{\rm)}   sphere $S^3(0,\epsilon)$, \, for $\epsilon > 0$
small enough. \item[{\rm(3)}] There are local analytic coordinates
$(x,y) \in (\bc^2,0)$ such that $\Omega(x,y)$ is either linear
$\Omega = \la xdy - \mu ydx$ with $\la/\mu \in \bc\backslash
\re_-$\,, or it is of the form $\Omega = xdy - (ny+x^n)dx$ where
$n \in \bn$ {\rm(}called Poincar\'e-Dulac normal form{\rm)}.
\end{enumerate}
\end{Theorem}

\noindent Thus we have the following problem:

\begin{Problem}
To obtain the local analytical  classification  of germs of
non-integrable holomorphic one-forms $\Omega$ at $0\in \bc^n$
under the hypothesis of transversality with small spheres
$S^{2n-1}(0,\epsilon)$.
\end{Problem}

One may work with the following notions and  models:

\begin{Definition} {\rm
 We shall say that a germ of singular non-integrable
one-form $\Omega$ at the origin $0 \in \bc^{2m}$ is {\it in the
Poincar\'e domain\/} if $\Ker(\Omega)$ is transverse to small
spheres $S^{4m-1}(0,\epsilon)$, \, $\epsilon
> 0$. A germ $\Omega$ will be called  {\it analytically linearizable\/}
if $f^*\Omega=\Omega_A$  for some germ of biholomorphism $f \in
\text{Bih}(\bc^{2m},0)$ fixing the origin and some $A\in \mathbb
A(2m)$. Finally we shall say that $\Omega$ is ({\it analytically
conjugate to}) {\it a non-integrable Poincar\'e-Dulac  normal
form\/} if $f^*\Omega = \Omega_{(\ell)}$  for some $f \in \text{
Bih}(\bc^{2m},0)$ and some $\ell \in \bn^m$.}
\end{Definition}

\begin{Question}
\label{Question:final}
 Are conditions {\rm (i)}, {\rm (ii)}  and {\rm (iii)}  below equivalent?

\noindent {\rm (i)}  $\Omega$ is the Poincar\'e-domain.

\noindent {\rm (ii)}  $\Omega$ is analytically linearizable or
conjugate to a non-integrable Poincar\'e-Dulac  normal form.

\noindent{\rm (iii)} There is a   holomorphic vector field $\xi$
transverse to the spheres $S^{4m-1}(0,r)$ for $r>0$ small enough
and such that $\Omega \cdot \xi=0$?
\end{Question}

\begin{Remark}
\label{Remark:final} {\rm The non-integrable examples we have
given, do admit holomorphic sections $\xi$ transverse to small
spheres $S^{4m-1}(0,r)$ (see Example~\ref{Example:distributions}).
Let us observe that there are examples of non-integrable one-forms
$\Omega$ with isolated singularity at the origin in $\bc ^n, \, n
\ge 3$, which admit integral manifold. For example, take
$f=\frac{1}{2}\sum\limits_{j=1}^{2m} z_j^2$ and $\Omega=df + f\nu$
for some holomorphic one-form $\nu=\sum\limits_{j=1}^{2m}\nu_j\,
dz_j$ such that $d\nu(0)$ is nondegenerate. Then $\{f=0\}$ is an
integral manifold of $\Omega$ and also $\Omega(z)=0$ if and only
if $z_j+ \sum\limits_{k=1}^{2m} z_k^2 \nu_j(z)=0$ for every $j$.
This shows, because the Jacobian of the left hand side at $z=0$ is
the identity, that $\Omega$ has an isolated singularity at the
origin $z=0$. Finally, if $X_j\in \Ker(\Omega)$, then
$(d\Omega)(X_1,X_2)=f(d\nu)(X_1,X_2)$, we obtain that near the
origin $\Omega$ is non-integrable at every point where $f\ne 0$.
This example suggests that a more interesting question might be :
\begin{Question}
\label{Question:classify} Classify the non-integrable germs of
holomorphic one-forms $\Omega$ with an isolated singularity at the
origin and $\Ker(\Omega)$ transverse to all small spheres centered
at the origin, which admit an integral manifold through the
origin.
\end{Question}

The above construction suggests that if the integral manifold is
taken in the form $\{f=0\}$ for a holomorphic function $f$, then
maybe we can write $\Omega=g(df + f\nu)$ for some function $g$ and
a one-form $\nu$, and we have to study the kernels of $df + f\nu$.
}
\end{Remark}

\subsection{Proof of Theorem~\ref{Theorem:isotopy}}

Now we shall  prove Theorem~\ref{Theorem:isotopy}. We need:
\begin{Lemma} \label{Lemma:arcwise}   The open subset $\ba(2m)\subset \mathcal A(2m)$
is arcwise-connected.
\end{Lemma}

\noindent{\bf Proof}.  First we recall that $\mathcal A(2m)$ is a
 complex vector space. Denote by $f\colon \mathcal A(2m) \to \bc$
the determinant function, i.e., $f(A) = \Det(A)$, $\forall\,A \in
\mathcal A(2m)$. Then $f$ is holomorphic and $f \not= 0$ because
if $\bj(2m) = \left[\begin{smallmatrix}
\big[\begin{smallmatrix}0 &-1\\ 1 &0\end{smallmatrix}\big]_{\ddots} &\\
&\big[\begin{smallmatrix}0 &-1\\ 1
&0\end{smallmatrix}\big]\end{smallmatrix}\right]$ then $\bj(2m)
\in \mathcal A(2m)$ and $f(\bj(2m)) = \Det(\bj(2m)) = (-1)^m \ne
0$.

\noindent Now $\ba(2m) = \big\{A \in \mathcal A(2m) : f(A) \ne
0\big\} = \mathcal A(2m)\backslash f^{-1}(0)$. Since $\mathcal
A(2m)$ is an affine space and $f^{-1}(0)$ is a thin subset of
$\mathcal A(2m)$ it follows from \cite{[Gunning-Rossi]} Corollary
4 page 20 that $\ba(2m)$ is (arcwise) connected.\qed

\vglue.1in We shall adopt the following natural definition:

\begin{Definition} \label{Definition 8.4} {\rm Two codimension one holomorphic
distributions $\Delta_1$, $\Delta_2$ on a complex manifold $M^n$
are {\it homotopic\/} if there is a $C^\infty$ family
$\{P_s\}_{s\in [0,1]}$ of holomorphic distributions in $M$ such
that $P_0=\Delta_1$ and $P_1=\Delta_2$. }
\end{Definition}

\begin{Lemma} \label{Lemma:linearhomotopy}  Let $\Omega_A$ be a linear one-form in
$\bc^{2m}$ with $A \in \ba(2m)$. Then the distribution
$\Ker(\Omega_A)$ is homotopic to the {canonical distribution}
$\Omega_{\bj(2m)} = \sum\limits_{j=1}^{2m} \big(z_{2j} dz_{2j-1} -
z_{2j-1} dz_{2j}\big)$.
\end{Lemma}

\noindent{\bf Proof}.  This is an immediate consequence of
Lemma~\ref{Lemma:arcwise}. Given a smooth path $\al\colon [0,1]
\to \ba(2m)$ connecting $\al(0) = A$ to $\al(1) = \bj(2m)$ we
define a homotopy  by setting $P_s=\Ker(\Omega_{\alpha(s)})$
(recall that $\al(s)\in \mathbb A(2m)$). \qed

\vglue.2in \noindent Now we consider a holomorphic one-form $\Om$
defined in a neighborhood $U$ of the closed unit ball $B^{4m}[0,1]
\subset \bc^{2m}$ and such that $\Om \cdot\vec R = 0$.  Write $\Om
= \sum\limits_{\nu=1}^{+\infty} \omega_\nu$ where $\omega_\nu$ is
a homogeneous one-form of degree $\nu \ge 1$ in $U$ (assume $U$ is
a ball). Then $0 = \Om\cdot\vec R = \sum\limits_{\nu=1}^{+\infty}
\omega_\nu\cdot\vec R$ and since $\omega_\nu\cdot \vec R$ is a
homogeneous polynomial of degree $\nu+1$ we conclude that
$\omega_\nu\cdot\vec R = 0, \quad \forall\,\nu \ge 1.$  In
particular, $\omega_1\cdot\vec R = 0$ and therefore $\omega_1 =
\Omega_A$ for some $A \in \mathcal A(2m)$. Assume now that
$\Ker(\Om)$ is transverse to $S^{4m-1}(1)$ then it follows from
Theorem~\ref{Theorem 8.1} that $\Omega$ has only one singularity
which is simple. Because the group of holomorphic transformations
of the unit ball acts transitively, we can assume by a holomorphic
change of coordinates that the origin is the only singularity of
$\Om$ in $B^{4m}[0,1]$ and, since it is a simple singularity, $A$
is non-singular, i.e., $A \in \ba(2m)$.

\begin{Lemma} \label{Lemma:deformation}  Under the above hypothesis there is a real
analytic deformation $\big\{\Ker(\Om^t)\big\}_{t\in[0,1]}$ of
$\Ker(\Om)$ into $\Ker(\Omega_A)$ by holomorphic distributions
$\Ker(\Omega^t)$ transverse to $S^{4m-1}(1)$ outside the
intersection $S^{4m-1}(1) \cap \Sing(\Omega^t)$.
\end{Lemma}

\noindent{\bf Proof}.  We define $\Om^t$ by $\Om^t
:=t^{-1}\Omega(tz)$. Then $\Om^t$ converges to a holomorphic
one-form for each $t \in [0,1]$ by canonical convergence criteria,
indeed, $\Om^t$ defines a holomorphic one-form in $U$ for each $t$
in the closed unit disc $\ov\bd \subset \bc$.  We have  $\Om^1 =
\Om$ and an easy computation shows that $\Omega^t$ converges for
$t\rightarrow 0$ to $\Omega_A$ where $A$ is the Jacobi-matrix of
$\Omega$ at the origin. Furthermore, because $\vec R(tz)=t\vec R
(z)$, one has $\Omega^t(z)\cdot \vec R(z)=t^{-1}\Omega(tz) \cdot
\vec R(z)= t^{-2} \Omega(tz) \cdot \vec R(tz) =0$ for every $t\ne
0$. Thus $\Ker(\Om^t)$ is transverse to $S^{4m-1}(1)$ outside
$\Sing(\Om^t) \cap S^{4m-1}(1)$.\qed

\vglue.1in

\noindent Now we study the singular set $\Sing(\Om^t)$.
\begin{Lemma}
\label{Lemma:singularset} We have $\sing(\Om^t) = \{0\}$,
$\forall\, t \in \ov\bd$.
\end{Lemma}
\noindent{\bf Proof}.  Take any $z\in S^{4m-1}(1)$. Consider the
complex line $\ell(z)=\{tz : t \in \bc\}$ through $z$. If $0<
|t|<1$, then $tz \ne 0$ belongs to $B^{4m}[1]$. Therefore
$\Omega^t(z)= \frac{1}{t}\Omega(tz)$ is non-zero. If $t=0$, then
$\Omega^0(z)=\Omega_A(z)$ is non-zero on $S^{4m-1}(1)$. This
proves the lemma. \qed

\vglue.1in \noindent{\bf Proof of Theorem~\ref{Theorem:isotopy}}:
For the proof we just have to apply
Lemmas~\ref{Lemma:linearhomotopy}, ~\ref{Lemma:deformation} and
~\ref{Lemma:singularset} above.\qed

\bibliographystyle{amsalpha}

\vglue.1in

\begin{tabular}{ll}
Toshikazu Ito  & \qquad  Bruno Sc\'ardua\\
Department of Natural Science  & \qquad Instituto de  Matem\'atica\\
Ryukoku University  & \qquad Universidade Federal do Rio de Janeiro\\
Fushimi-ku, Kyoto 612 & \qquad  Caixa Postal 68530\\
JAPAN   & \qquad 21.945-970 Rio de Janeiro-RJ\\
&  \qquad BRAZIL
\end{tabular}

\end{document}